\documentclass[11pt]{article}
\usepackage{latexsym}
\usepackage{amssymb, amsmath }
\usepackage[dvips]{graphicx}
\usepackage{amsfonts}
\usepackage{fancyhdr}
\usepackage{amscd}
\newtheorem{teo}{Theorem}[section]
\newtheorem{prop}[teo]{Proposition}
\newtheorem{lema}[teo]{Lemma}
\newtheorem{coro}[teo]{Corollary}
\def\C{{\mathbb C}}
\def\N{{\mathbb N}}
\def\Q{{\mathbb Q}}
\def\R{{\mathbb R}}
\def\Z{{\mathbb Z}}

\newcommand{\ds}{\displaystyle}

\newcommand{\diam}{\operatorname{diam}}

\newcommand{\n}{\hspace{0pt}}
\newcommand{\la}[1]{\lambda_{#1}}
\newcommand{\La}{\Lambda}

\newcommand{\p}{P_{\lambda}}
\newcommand{\q}{Q_{\lambda}}
\newcommand{\D}{\textbf{Proof:}}

\begin{document}

\title{WANDERING FATOU COMPONENTS ON $p$-ADIC POLYNOMIAL DYNAMICS}
\author{Gabriela Fern\'andez Lamilla}
\date{}

\maketitle \thispagestyle{empty} \abstract{}
 We will study perturbations of the polynomials
$P_\lambda$, of the form $$Q_\lambda= P_\lambda +Q$$ in the space
of  centered monic polynomials, where $P_\lambda$ is the
polynomial family defined by
$$P_{\lambda}(z)=\frac{\lambda}{p}z^p+\left(1-\frac{\lambda}{p}\right)
z ^{p+1}$$ with $\lambda \in \Lambda= \{z: |z-1| <1\}$, studied
by Benedetto , who showed that for a dense set of parameters, the
polynomials  $P_\lambda$ have a wandering disc contained in the
filled Julia set.

We will show an analogous result for the family $Q_\lambda$,
obtaining the following consequence:

The polynomials $P_\lambda$ belong to
$\overline{\mathrm{E}}_{p+1}$ where $\mathrm{E}_{p+1}$ denotes
the set of polynomials that have a wandering disc in the filled
Julia set, in the space of centered monic polynomials of degree
$p+1$.
\newpage
\thispagestyle{empty}
\tableofcontents
\newpage

\setcounter{page}{1}
\section{Introduction}

\indent \indent In complex dynamics there exists an extensive
study of polynomials as dynamical systems acting on $\C$. The
orbit of a point $z_0\in\C$ under a polynomial $f \in \C[z]$ is
the sequence $z_0,z_1,z_2,\ldots$ defined by
$$z_n =f^n(z_0).$$

\n Subsets of $\C$ which are of particular interest are the
\textbf{filled Julia set}, which is the set of points with bounded
orbit; the \textbf{Julia set}, which corresponds to the boundary
of the filled Julia set; and the \textbf{Fatou set}, the
complement of the Julia set.

\n A important result of Sullivan \cite{Su} says that for all
polynomials $f \in \C[z]$ there is no wandering component of the
Fatou set, i.e. every connected component of the Fatou set is
pre-periodic under the action of $f$ (This result holds also for
rational functions but our emphasis will be on polynomials).

\n Recently the study of iterations of rational functions over
$\C$ has been extended to the study of rational functions with
coefficient in the field $\C_p$ (\cite{BE}, \cite{BE3},
\cite{RL2}, \cite {RL3}). This field is the smallest complete
algebraically closed extension of $\Q$ with respect to the
$p$-adic valuation. The construction of $\C_p$ is analogous to
that of the complex numbers starting with rational numbers and the
usual absolute value, some interesting differences arise between
$\C$ and $\C_p$.

\n The field $\C_p$, endowed with the $p$-adic valuation, is an
\textbf{ ultrametric space}, i.e. for all $x,y \in \C_p$
$$|x+y|\leq \max \{|x|,|y|\}.$$

\n From the above inequality, known as the \textbf{strong triangle
inequality}, it follows that $\C_p$ is totally disconnected, so
the connected component notion used in complex dynamics must be
replaced by the concept of \textbf{infraconnected component} (see
\cite{ES}).

\n The motivation of this work arises from a result of Benedetto
\cite{BE}, who studied the family of polynomials in $\C_p[z]$
defined by
$$P_{\lambda}(z)= \frac{\lambda}{p} \ z^p
+\left(1-\frac{\lambda}{p}\right)\ z^{p+1},$$ where $\lambda \in
\Lambda =\{\lambda \in \C_p: |\lambda-1|_p<1\}$, obtaining the
following result:

\bigskip
\noindent {\bf Theorem (Benedetto).} {\em There is a dense set of
parameters $\lambda \in \Lambda$ such that the polynomial
$P_\lambda $ has a wandering disc contained in the filled Julia
set which is not attracted to an attracting cycle.}
\bigskip

\n From the above theorem we conclude that there exist polynomials
in $\C_p[z]$ with wandering infraconnected components of the Fatou
set, in contrast with the result of Sullivan for complex rational
functions.

\n In this work, we will study  perturbations of the polynomials
$P_\lambda$, of the form $$Q_\lambda = P_\lambda + Q,$$ where $Q$
is a polynomial with $\|Q\|\leq
\left(\frac{1}{p}\right)^{\frac{1}{p-1}}$ and we will obtain the
following result (for the definition of $\|Q\|$ see Section 2.1).

\bigskip
\noindent {\bf Theorem.} {\em There is a dense set of parameters
$\lambda \in \Lambda$ such that $Q_\lambda $ has a wandering disc
contained in the filled Julia set which is not attracted to an
attracting cycle.}
\bigskip

\n Let $Pol_d$ be the space of monic centered polynomials of
degree $d\geq 2$ and coefficients in $\C_p$. The parameter space
$Pol_d$  is naturally identified with $\C_p^{d-1}$. If we call
$E_d$ the set of polynomials in $Pol_d$ that have a wandering
disc, from the above theorem, we obtain directly the following
consequence.

\bigskip
\noindent {\bf Corollary.} {\em For all $\lambda \in \Lambda$, the
polynomial $ P_\lambda $ belong to the interior of
$\overline{E}_{p+1}$.}
\bigskip

\n In addition, we will prove that the above theorem is also true
for a wider class of perturbations of the polynomials $\p$. In
fact, if we consider $R_B$ the set of rational functions without
poles in the fixed ball $B= \{z: |z|\leq r\}\ \ (r>1)$ and the
subset $R_B^{E}$ of functions in $R_B$ with a wandering disc, we
obtain the following consequence.

\bigskip
\noindent {\bf Corollary.} {\em For all $\lambda \in \Lambda$, the
polynomial $ P_\lambda $ belong to the interior of
$\overline{R_{B}^E}$.}

\medskip
\n In sections 2.1 and 2.2 we recall some basic concepts and facts
from ultrametric analysis and dynamics. In Section 2.3 we will
present in detail some results and techniques used in \cite{RL}
since they are essential for our study of the perturbation
$Q_{\lambda}$. Our study is not done directly on $\q$, but it is
more convenient to work with an affinely conjugated map
$Q_{\lambda}^*$. In Section 3.1 we study the behavior of
$Q_{\lambda}^*$ over the filled Julia set. Finally, in Section 3.2
we mimic \cite{BE} to study $Q_\lambda^*$ as a function of the
parameter $\lambda$ and prove our main result.

\newpage
\section{Preliminaries}

\indent \indent In this section we recall definitions and results
that are used throughout this work.

\n The field $\C_p $ endowed with the $p$-adic valuation denoted
by $|\cdot|$ is an ultrametric space, i.e. for $z_0,z_1 \in \C_p$
we have that $|z_0-z_1|\leq\max \{|z_0|,|z_1|\}$. From this
inequality and the completeness of $\C_p$ arise interesting
topological and geometrical results, some of them are:
\begin{itemize}
\item[(i)] The {\bf value group} of the valuation is the set
$|\C_p^*|:=\{|z|:z\in \C_p^{*} \}=\{p^r:r\in\Q \}$.

\item[(ii)]{\bf The isosceles triangle principle}:  If
$|z_1|\neq|z_2|$, then $|z_1+z_2|=\max \{ |z_1|,|z_2|\}$.

\item[(iii)] For $z_0 \in \C_p$, $r\in |\C_p^{*}|$, the {\bf open
ball} with radius $r$ and center $z_0$ is the set
$$B_r(z_0)=\{ z \in \mathbb{C}_p : |z-z_0|< r\}$$ and the {\bf closed ball}
with center $z_0$ and radius $r$ is the set
$$\overline{B_{r}}(z_0)=\{ z \in \mathbb{C}_p: |z-z_0|\leq r \}.$$ These are open and closed sets in the topology of
$\mathbb{C}_p$. By definition, we have that the diameter
 of $B$ denoted $\diam(B)$ belong to $|\C_p^*|$, where $B$ is an open or closed ball.

\n We denote by $\mathcal{O}_{\C_p}$ the closed ball $
\{z:|z|\leq1\}.$

\item[(iv)] Every point of a ball is a center, that is, if  $z_1
\in B_r(z_0)$ (resp. $z_1 \in \overline{B_r}(z_0)$), then
$B_r(z_1) =B_r(z_0)$ (resp.
$\overline{B_r}(z_1)=\overline{B_r}(z_0)$).

\item[(v)] If two balls have not empty intersection then one is
contained in the other one.
\end{itemize}

\n The properties below are about convergence in ultrametric
spaces, some of them are different than in archimedean analysis.
Let $a_1, a_2,\dots,a_n,\dots $ be a sequence in $\C_p$. Then
\begin{itemize}
\item[(vi)] If $\underset{n\rightarrow \infty}{\lim} a_n=a$ and
$a\neq 0$, then there exists a $n_0$ such that $|a_n|=|a|$ for all
$n>n_0$.

\item[(vii)]The series $\underset{n=1}{\overset{\infty}{\sum}}a_n$
converges if and only if $\underset{n\rightarrow \infty}{\lim} a_n
=0$.

\item[(viii)] The power series
$\underset{n=1}{\overset{\infty}{\sum}}a_nx^n$ has convergence
radius $r:= \left(\underset{n \rightarrow \infty}{\limsup}
\sqrt[n]{|a_n|}\right)^{-1}$.

\item[(ix)] If $r$ is the convergence radius of
$\underset{n=1}{\overset{\infty}{\sum}}a_nx^n$, then the map

$$x\longmapsto
\underset{n=1}{\overset{\infty}{\sum}}a_nx^n$$

\n is  differentiable in $B_r(0)$ and its derivative is
$$x\longmapsto
\underset{n=1}{\overset{\infty}{\sum}}na_nx^{n-1}.$$
\end{itemize}

\subsection{Ultrametric analysis.}
\indent \indent Let $B$ be a ball with radius $r$ (open or
closed), we denote by $\mathcal{H}(B)$ the ring of power series
which converge in $B$. The space $\mathcal{H}(B)$ endowed with the
norm
$$\| f\|_{B}=\underset{i\geq0}{\sup}|a_i|r^{i}$$  is a complete ultrametric valued ring.

\n As in the complex case, a rational function can be written as a
power series around every point $z_0$ which is not a pole. But
observe that if we consider the function $f:\C_p \longrightarrow
\C_p$ given by
$$f(z)=\left\{ \begin{array}{ll}
   1, & |z|<1\\
   0, & |z|\geq 1 \\
\end{array} \right.$$ we have that also can be written as power series around every point of $\C_p$.
Hence, it is clear that the idea of holomorphic functions in
$\C_p$ is different from the one in $\C$, which it is defined by a
local property. Indeed, a function defined in a subset $X$ of
$\C_p$ is holomorphic if it is the uniform limit of rational
functions without poles in $X$ (see \cite{TJ}). In this work, we
only consider holomorphic functions defined on a ball $B$ and, in
this case, the definition coincides with the complex one: a
function $f$ is \textbf{holomorphic} in ${B}$ if and only if $f$
can be written as a convergent power series in $B$. Thus,
$\mathcal{H}(B)$ is the space of holomorphic functions in the disc
$B$.

\bigskip
\n Now, we will show an analogous to the Newton's method, in order
to guarantee the existence of roots in a holomorphic function.

\begin{lema}{\bf(Hensel)} {\hspace*{-3pt}\bf .} Let $f\in \mathcal{O}_{\C_{p}}[[z]]$.
If there exists $z_{0}\in \mathcal{O}_{\C_p}$ with $|f(z_{0})|<
|f^{'}(z_0)|^2$, then there is an unique root  $w$ of $f$ such
that $|w-z_0|\leq \ds
\frac{|f(z_0)|}{|f^{'}(z_0)|}$.\label{Hensel}
\end{lema}
\D

\n We define recursively the sequence,
$$z_1= z_0 -\frac{f(z_0)}{f^{'}(z_0)},
\ z_{n+1}= z_n-\frac{f(z_n)}{f^{'}(z_n)}.$$ We will show
inductively that:
\begin{itemize}
\item[(i)] $|f(z_n)|\leq C^{2^{n}}|f^{'}(z_0)|^{2}$, where $C=
\ds\frac{\left |f(z_0)\right |}{\left |f^{'}(z_0)\right |^{2}}<1.$
\item[(ii)]$|f^{'}(z_n)|=|f^{'}(z_0)|.$
\item[(iii)]$|z_n-z_0|=|z_1-z_0|.$
\end{itemize}

\n Now
$$f(z_1)=f(z_0+z_1-z_0)=f(z_0)+f'(z_0)(z_1-z_0)+d (z_1-z_0)^{2}=
  d(z_1-z_0)^2$$ for some $d=d(z_0,z_1) \in
  \mathcal{O}_{\C_p}$. Therefore

\begin{equation*}
|f(z_1)| = |d(z_1-z_0)^2| \leq  |z_1-z_0|^2 = C^2|f'(z_0)|^2.
\end{equation*}
and
 $$f'(z_1)= f'(z_0)+e (z_1-z_0)$$ for some $e=e(z_0,z_1)\in \mathcal{O}_{\C_p}$. Hence $$|f'(z_1)-f'(z_0)|\leq |z_1-z_0|
  \leq\left | \ds \frac{f(z_0)}{f'(z_0)}\right |<|f'(z_0)|.$$

From the previous inequality, and the isosceles triangle principle
applied to  $|f'(z_1)-f'(z_0)+f'(z_0)|$ we have that
$$|f'(z_1)|=|f'(z_0)|.$$

\n The inductive steps for (i) and (ii) are analogous to the
previous one, (iii) is direct consequence  of (i) and (ii).

\n Since $|f'(z_n)|=|f'(z_0)|$, we have that $|z_{n+1}-z_n|=\ds
\frac{|f(z_n)|}{|f'(z_n)|}
  \leq C^{2^n}|f'(z_0)|$, so $\{ z_n\}_{n\in \mathbb{N}}$ is a Cauchy sequence and if we
denote its limit by $w$ we get, from (i) and (iii), that $w$ is a
root of
 $f$ and $|w-z_0|\leq \ds \frac{|f(z_0)|}{|f^{'}(z_0)|}$. \hfill $\Box$

\bigskip

\n We now enumerate some interesting properties of holomorphic
functions (see \cite{BE2}).
\begin{teo}{\hspace*{-6pt}\bf .}
Let $f(z)=\overset{\infty}{\underset{i=0}{\sum}}a_iz^{i}$ with
$a_i \in \mathbb{C}_p$ and $r \in |\mathbb{C}_p^{*}|$ such that
$\underset{i\rightarrow\infty}{\lim}|a_i|r^{i}=0$. Then $f$ has a
root $\alpha \in \mathbb{C}_p$ with $|\alpha|=r$ if only if there
exist $n, m \in \Z $ with $n<m$ and such that
\begin{equation}
|a_n|r^n=|a_m|r^m=\underset {i\geq 0}{\sup}\{|a_i|r^{i}
\}\label{ec m y n}.
\end{equation}

\n Moreover if $n,m$ are the smallest and the greatest integers,
respectively, that make \emph{(\ref{ec m y n})} true, then $f$ has
exactly $m-n$ roots with absolute value $r$, counting
multiplicity.\label{teorema de soluciones contar}
\end{teo}

\begin{coro}{\hspace*{-6pt}\bf .}
Let $B$ be a closed ball, $f\in \mathcal{H}(B)$, and $D$ an open
ball (resp. closed) contained in $B$. Then $f(D)$ is an open ball
(resp. closed).\label{imagen de bola}
\end{coro}

\begin{coro}{\hspace*{-6pt}\bf .}
Let $f\in \mathcal{H}(B)$, $w_0 \in \mathbb{C}_p$, $\delta \in
|\mathbb{C}_p^{*}|$ such that $\overline{B_{\delta}}(w_0)\subset B
$. If $|f(w)-f(w_0)|= \alpha$ for all $w$ in $\{ w:|w-w_0|=\delta
\}$, then
$$f(\{ w:|w-w_0|=\delta\})=\{w:|w-f(w_0)|=\alpha \}.$$\label{esfera}
\end{coro}

\begin{coro}{\hspace*{-6pt}\bf .} Let $B$ be a closed ball and $f\in \mathcal{H}(B)$.
Then there exists $d\in\N$ such that the series $f- w$ has exactly
$d$ roots in $B$, counting multiplicity, for all $w \in f(B)$.
\end{coro}

\n For $f\in \mathcal{H}(B)$, we define the \textbf{degree of the
map} $f$ as the number $d$ from the last corollary.

\subsection{Polynomial dynamics over $\C_p.$}
\indent\indent This section contains some important definitions
and dynamical properties that will be needed later.

\bigskip
\n Let $P \in \C_p[z]$. For $z \in \C_p$, we define the
\textbf{orbit} of $z$, denoted by $\mathcal{O}(z)$, as the
sequence $\{P^n(z)\}_{n\in \N}$. If $P(z)=z$, we say that $z$ is a
\textbf{fixed point} of $P$; if for $z$ there exists $n\in \N$
such that $P^n(z)=z$, we will say $z$ is a \textbf{periodic
point}.

\n Let $P'$ be the formal derivative of $P$, $z_0$ a fixed point
of $P$ and $\theta =|P'(z_0)|$. Then:
\begin{itemize}
\item[(i)] If $\theta < 1$, we say that $z_0$ is an
\textbf{attracting fixed point}. \item[(ii)] If $\theta > 1$, we
say that $z_0$ is a \textbf{repelling fixed point}. \item[(iii)]
If $\theta =1$, we say that $z_0$ is an \textbf{indifferent fixed
point}.
\end{itemize}

\bigskip
\n Other object of study is the \textbf{filled Julia set} denoted
by $K(P)$, that correspond to the set of points of $\C_p$ with
bounded orbit. Some properties of the filled Julia set are:
\begin{enumerate}
\item[(i')] $K(P) \neq \emptyset$. \item[(ii')] $K(P)$ is closed
and bounded. \item[(iii')]$P^{-1}(K(P))= K(P)$, i.e. $K(P)$ is
completely invariant.
\end{enumerate}
\bigskip

\n Another important set is the \textbf{Julia set}, $J(P)$, which
is the boundary of the filled Julia set. The Julia set can be also
defined as follows
$$\{z \in \C_p: \text{ for every neighbourhood } U \text{ of } z,\,
\underset{n \in \N}{\bigcup}P^n (U)= \C_p\}.$$

\n Finally, we define the \textbf{Fatou set} as the complement of
the Julia set. We denote it by $F(P)$.

\bigskip
\n We are not going to study just polynomials, so we have to
introduce the concept of polynomial like maps. If $U$ and $V$ are
open balls in $\C_p$ such that $U\subsetneq V$ and
$f:U\longrightarrow V$ is a holomorphic function of degree $d$
with $d\geq 1$, we say that $(f, U)$ is a \textbf{polynomial like
map} of degree $d$.

\n All the preceding concepts can be also defined, in a similar
way, for polynomial like maps. That is, the filled Julia set of
$(f,U)$ is the set
$$K(f,U)=\{z\in U: f^n(z)\in U \text{ for all } n\in \N\}.$$
\n The Julia set is $$J(f,U)=\partial K(f,U),$$ and the Fatou set
is
$$F(f,U)= U\setminus J(f,U).$$

\n These new definitions will allow us to study the dynamical
behavior of some holomorphic functions restricted to balls.

\n Let $D$ be a subset of $\C_p, \,a \in D$ and $I_a:\C_p
\longrightarrow \R$ the map defined by $I_a(x)=|x-a|$. We say that
$D$ is \textbf{ infraconnected} if and only if for all $a \in
\C_p$ the set $\overline{I_a(D)}$ is an interval (see \cite{ES}).
In particular, for $(f,U)$, we are interested in understanding the
behavior of the filled Julia set. If we consider $B$, the smallest
ball that contains $K(f,U)$, then $f^{-n}(B)$ is a collection of
disjoint closed balls, named balls of level $n$. Then for $w\in
K(f,U)$ there is an unique sequence $\{B_n\}_{n\in\N}$ of nested
closed balls, where $B_n$ is a ball of level $n$, such that $w\in
\underset{n\in \N}{\cap} B_n$. The set $C(w):=\underset{n}{\cap}
B_n$ is the infraconnected component of $K(f,U)$ that contains $w$
(see \cite{TJ}).

\n Now let $(f,U)$ be a polynomial like map. We say that
$E\subseteq U$ is a \textbf{wandering set} if $f^n(E)\bigcap
f^m(E)\neq \emptyset$ only when $n=m$.

\n Furthermore, if $(f,U)$ and $(g,U)$ are polynomial like maps,
if there is a homomorphism $h:U\longrightarrow U$ such that
$g=h^{-1}f h$, we say that $f$ and $g$ are \textbf{topologically
conjugated}. In these case, if $E$ is a wandering set of $f$, then
$h(E)$ is a wandering set of $g$. Therefore, the existence of
wandering set is invariant under conjugacy. This fact will turn
out to be very important to obtain our results.

\newpage
\subsection{The family of polynomials $P_\lambda$.}

\indent\indent For $\lambda \in \Lambda = \{ \lambda \in
\mathbb{C}_p : |\lambda -1|<1\}$ let
$$P_{\lambda}(z)= \frac{\lambda}{p}z^p
+\left( 1-\frac{\lambda}{p} \right) z^{p+1}.$$
\begin{teo}[Benedetto]{\hspace*{-6pt}\bf .} There is a dense set of parameters $\lambda \in \Lambda$, such that
the polynomial $P_{\lambda}$ has a wandering disc contained in
$K(P_{\lambda})$, which is not attracted to an attracting cycle.
\end{teo}

\n Now we will sketch the proof of this theorem (see \cite{RL}),
paying attention to the techniques which will be important later.

\n First we notice that $\overline{B_{\rho}}(0)$, with $\rho = \ds
p ^{\frac{-1}{p-1}}$, is invariant under the action of $\p$ and
that $z=1$ is a repelling fixed point. Now, we considere $B_1(0)$
and $B_1(1)$, which are neighbourhoods of the fixed ball
$\overline{B_\rho(0)}$ and the repelling fixed point $z=1$
respectively. From the strong triangle inequality and Corollary
\ref{imagen de bola}, we see that the set $K(\p)$ is contained in
$B_1(0) \sqcup B_1(1)$. This allows us to define the itinerary of
a point $x \in K(\p )$ as the sequence
$$\theta _1 \theta_2 \ldots \theta _n \ldots $$ with
$\theta_i \in \{0,1\}$ and $  \p^{i}(x) \in B_1(\theta_i)$ for $i
\in \N$. Furthermore, we obtain that all the points of a ball
contained in the filled Julia set have the same itinerary. If this
itinerary is not pre-periodic, then $D$ is a wandering disc. In
order to find such disc is necessary to study the behavior of the
$\p$ in the filled Julia set. The lemmas below describe such
behavior.

\n  We define $S>0$ by $pS^{p-1}=\rho$ and the sequence $\{\rho
_n\}_{n \in \N}$ by
$$\rho_{0} =1,\ \ p\rho_{n}^p = \rho_{n-1}.$$  \vspace*{-60pt}
\begin{lema}{\hspace*{-6pt}\bf .}
$\begin{array}{rl} \\*[55pt] 1)& \text{ Let } m\geq 1, z_0 \text{
and } z_1 \text{ such that }|z_{0}|=| z_{1}|=\rho _{m}.\text{ If }
| z_{0}-z_{1}| \leq S, \text{ then: }\\*[8 pt] & \hspace*{2
cm}|P_{\lambda}(z_{0})-P_{\lambda}(z_{1})| \leq \rho _{m-1}
|z_{0}-z_{1}|.\\*[8 pt] 2) & \text{ If }  z_{0},z_{1} \in
B_{1}(1), \text{ then: } \\*[8pt] &\hspace*{2
cm}|P_{\lambda}(z_{0}) -P_{\lambda}(z_{1})| =p\ |z_{0}-z_{1} |.
\end{array}$\label{comp de P}
\end{lema}
 \D

\begin{itemize}
\item[{\em1)}]We observe that

$P_{\lambda}(z_{0})-P_{\lambda}(z_{1}) =\ds \frac{\lambda}{p}\left
(p\, \varepsilon \,z_0^{p-1}+\cdots
+p\,\varepsilon^{p-1}z_0+\varepsilon^p\right
)+\left(1-\frac{\lambda}{p}\right )((p+1)\,\varepsilon\,
z_0^p+\cdots +\varepsilon ^{p+1}),$

with $\varepsilon =z_1-z_0$, since $|\varepsilon |= |z_0-z_1|\leq
S< \rho_m$ and

$$|P_{\lambda}(z_0)-P_{\lambda}(z_1)|\leq |\varepsilon |\max
\{z_0^{p-1}, p\,|\varepsilon |^{p-1},\rho_{m-1}\}= |\varepsilon |\
\rho_{m-1},$$

we have that $$|P_{\lambda}(z_0)-P_{\lambda}(z_1)|\leq \rho_{m-1}\
|z_0-z_1|.$$

 \item[{\em2)}] The proof is straightforward from the previous one and will be omitted.\hfill
$\Box$
\end{itemize}

\n With the previous lemma it is possible to find a necessary and
sufficient condition for the existence of wandering discs in
$K(\p)$, this condition is:

\begin{lema}{\hspace*{-6pt}\bf .}
Let $\{ m_{i}\} _{i\geq 0},\{ M_{i}\} _{i\geq 0}$ be two sequences
of positive integers such that, for all $i\geq0$ we have that
$\rho_{m_{i}-1} \cdot \ldots \cdot \rho_{1}\cdot p^{M_{i}}\leq1.$
Suppose that for $\lambda_{0} \in \Lambda $ there exists $x\in
K(P_{\lambda_{0}})$ with itinerary
$$\underset{m_0}{\underbrace{0\ldots0}}\,\underset{M_0}{\underbrace{1\ldots1}}\,
\underset{m_1}{\underbrace{0\ldots0}}
\,\underset{M_1}{\underbrace{1\ldots1}}\ldots\underset{m_i}{\underbrace{0\ldots0}}\,\underset{M_i}{\underbrace{1\ldots1}}\ldots,$$
then the ball $U=\{z:|z-x|\leq S \}$ is contained in
$K(P_{\lambda_{0}})$.\label{ex para P}
\end{lema}

\n Therefore, to prove Theorem 2.6 it suffices to find $x \in
K(\p)$ and sequences $\{M_i\}_{i \in \N},\{m_i\}_{i \in \N}$  with
$\lim M_i = \infty$ such that the hypothesis of the previous lemma
are satisfied. In order to do this we study the function $\p(z)$
as a function of $\lambda$.

\n Now, we will see two lemmas that will allow us to find such
sequences $\{M_i\}_{i\in \N},\{ m_i\}_{i\in \N}$ with $\lim M_i=
\infty$ implying the existence of wandering discs in $K(\p)$.

\begin{lema}{\hspace*{-6pt}\bf .}
Let $m\geq1$, and $z_0,z_1 \in B_{p^{-m}}(1)$. If $\lambda_0 ,
\lambda_1 \in \Lambda$ satisfy $|z_0 -z_1|=|\lambda_0
-\lambda_1|$, then \label{sn1}
$$|P_{\lambda_0}^{m}(z_0)-P_{\lambda_1}^{m}(z_1)|=p^{m}|\lambda_0 -\lambda_1|.$$
\end{lema}
\D  \ We proceed by induction. From
\begin{eqnarray}
&&|P_{\lambda_0}(z_0)-P_{\lambda_0}(z_1)|=p\,|z_0-z_1|=p\,|\lambda_0
-\lambda_1|, \label{eq1}\\*[.3 cm]
&&|P_{\lambda_0}(z_1)-P_{\lambda_1}(z_1)|=p\,|\lambda_0
-\lambda_1||z_1-1|< p\,|\lambda_0-\lambda_1|, \label{eq2}
\end{eqnarray}
we have that $|P_{\lambda_0}(z_0)-P_{\lambda_1}(z_1)|= p\,
|\lambda_0 -\lambda_1|$, therefore the lemma is true for $m=1$.

\n Now, for the inductive step, we suppose that $z_0, z_1 \in \{
z: |z-1|\leq p^{-m}\}$. By hypothesis we have that
$$|P^{m-1}_{\lambda_0}(z_0)-P^{m-1}_{\lambda_1}(z_1)|=
p^{m-1}|\lambda_0- \lambda_1|.$$
$$|P^{m-1}_{\lambda_1}(z_1)-1|<p^{m-1}.$$
\n Therefore,
\begin{eqnarray}
&&|P_{\lambda_0}(P^{m-1}_{\lambda_0}(z_0))-P_{\lambda_0}(P^{m-1}_{\lambda_1}(z_1))|=p^{m
}|\lambda_0-\lambda_1|,\label{eq3}\\*[.3cm]
&&|P_{\lambda_0}(P^{m-1}_{\lambda_1}(z_1))-P_{\lambda_1}(P^{m-1}_{\lambda_1}(z_1))|=
p\,|\lambda_0-\lambda_1||P_{\lambda_1}^{m-1}(z_1)-1|<p^{m}|\lambda_0
-\lambda_1|.\label{eq4}
\end{eqnarray}

\n From (\ref{eq3}) and (\ref{eq4}) we obtain
$$|P_{\lambda_0}^{m}(z_0)-P_{\lambda_1}^{m}(z_1)|=p^{m}|\lambda_0
-\lambda_1|.$$
\hfill $\Box$

\begin{lema}{\hspace*{-6pt}\bf .}
Let $m\geq 1 $ and $z_0, z_1$ with $|z_0|=|z_1|=\rho _m$ and such
that $|z_0-z_1|\leq S.$ If $\lambda_0 ,\lambda_1\in \Lambda$ are
such that
$$\rho_{m-1}\cdot \ldots\cdot \rho_1  \cdot |z_0-z_1|< |\lambda_0-
\lambda _1|\leq S,$$ then
$$|P_{\lambda_0}^{m}(z_0)-P_{\lambda_1}^{m}(z_1)|=|\lambda_0
-\lambda_1|.$$ \label{lema en P inductivo}
\end{lema}
\n \D \, First we will show inductively that, for $1\leq i\leq m$,

\begin{eqnarray}
&&|P_{\lambda_0}^{i}(z_0)-P_{\lambda_1}^{i}(z_1)|\leq \max \{
\rho_{m-1}\cdot\ldots\cdot \rho_{m-i}|z_0-z_1|,
\rho_{m-i}|\lambda_0-\lambda_1|\}.\label{eq6}
\end{eqnarray}
Observe that
\begin{eqnarray}
&&|P_{\lambda_0}(z_0)-P_{\lambda_0}(z_1)|\leq
\rho_{m-1}|z_0-z_1|,\label{eq7}
\\&&|P_{\lambda_0}(z_1)-P_{\lambda_1}(z_1)|=
\rho_{m-1}|\lambda_0-\lambda_1|.\label{eq8}
\end{eqnarray}
Using the ultrametric inequality, (\ref{eq7}) and (\ref{eq8}), we
get (\ref{eq6}) for $i=1$.

If we assume (\ref{eq6}) as the inductive hypothesis, we have
\begin{eqnarray}
&&|P_{\lambda_0}(P^{i}_{\lambda_0}(z_0))-P_{\lambda_0}(P^{i}_{\lambda_1}(z_1))|\leq
\rho_{m-i-1}|P^{i}_{\lambda_0}(z_0)-P^{i}_{\lambda_1}(z_1)|,\label{eq9}\\
&&|P_{\lambda_0}(P^{i}_{\lambda_1}(z_1))-P_{\lambda_1}(P^{i}_{\lambda_1}(z_1))|=\rho_{m-
i-1}|\lambda_0-\lambda_1|.\label{eq10}
\end{eqnarray}

\n From (\ref{eq9}) and (\ref{eq10}) we obtain (\ref{eq6}) for
$i+1$. Notice that for the inductive step from $m-1$ to $m$, the
hypothesis of the lemma gives us that
$$|P_{\lambda_0}^{m}(z_0)-P_{\lambda_1}^{m}(z_1)|=|\lambda_0
-\lambda_1|.$$ \hfill $\Box$

\n Let $\lambda \in \Lambda$ and $M_0 \in \N$ with $p^{-M_0} \leq
S$, we choose $m_0 \in \N$ such that $$\rho_{m_0-1} \cdot \ldots
\cdot \rho_1 p^{M_0} \leq 1.$$

\n Now, if we choose $x \in K(P_\lambda)$ with itinerary
$$\underset{m_0}{\underbrace{0\ldots 0}}\,1\,1\,1\,1\,1 \ldots,$$
we obtain Lemma \ref{sn1} hypothesis with $z_0 =
P_{\lambda_0}^{m_0}(x), z_1=P_{\lambda_1}^{m_0}(x)$ and $M=M_1$,
for all $\lambda_0, \lambda_1 \in \{z: |\lambda- z|\leq
p^{-M_0}\}$, and we have
$$|P_{\lambda_0}^{m_0+M_0}(x)-P_{\lambda_1}^{m_0+M_0}(x)|=
p^{M_0}|\lambda_0 - \lambda_1|.$$ Hence, there exists $w_0 \in
\Lambda$ with $P_{w_0}^{m_0+M_0}(x)=0$ such that the itinerary of
$x$ for $P_w$ is
$$\underset{m_0}{\underbrace{0\ldots 0}}\,\underset{M_0}{\underbrace{1\ldots 1}}\,0 \ldots$$
for all $w \in \{z: |z-w_0|< p^{-M_0}\}$.

\n As before, we choose $M_1$ such that $P^{M_1-M_0}\leq S$, and
$m_1$ such that
$$\rho_{m_1-1} \cdot \ldots
\cdot \rho_1 p^{M_1} \leq 1$$ obtaining that there exists
$\lambda' \in \Lambda$ with $|\lambda'- w_0|= \rho_{m_1}p^{-M_0}$
such that $P_{\lambda'}(x)=1$.

\n Therefore, the itinerary of $x$ for $P_{\lambda'}$ is
$$\underset{m_0}{\underbrace{0\ldots
0}}\,\underset{M_0}{\underbrace{1\ldots
1}}\,\underset{m_1}{\underbrace{0\ldots 0}}\,1\,1\,1 \ldots$$

\n Lemma \ref{lema en P inductivo} allows us to make this process
inductively, obtaining the Lemma \ref{ex para P} hypothesis.

\section{Results.}

\indent\indent In this section, we establish some properties of
the perturbations of the polynomials $P_\lambda$. Throughout,
$$\rho_0 =1,\ \ p\ \rho_n^p =\rho_{n-1}.$$

\n Recall that $ p\rho^p=\rho$ and that $pS^{p-1}=\rho$. For the
rest of this work we fix $\widehat{r}\in |\C_p^*|$, with
$\widehat{r}>1$ and $B=\{z\in\C_p:
 |z|\leq \widehat{r}\}$. The perturbations are:

$$Q_\lambda^*(z)=P_\lambda(z)+Q(z),$$
 where $Q\in \mathcal{H}(B)$ with $\|Q\|< \rho$. For this family we will obtain the following result:
\begin{teo}{\hspace*{-6pt}\bf .}
There is a dense set of parameter $\lambda \in \Lambda$ such that
the function $Q_{\lambda}^*$ has a wandering disc contained in the
filled Julia set, which is not attracted to an attracting
cycle.\label{teo pa estrella}
\end{teo}

\n To prove this theorem, we will study a topological conjugation
of $Q_\lambda^*$.

 \n Notice that $$p(Q_\lambda^*(z)-z) \in \mathcal{O}_{\C_p}[[z]],$$
 in addition $$\left|p\ Q_\lambda^*(1)-p \,\right|=\frac{1}{p}\
 |Q(1)|<\frac{\rho}{p}$$ and
 $$|p\ (Q^*_\lambda)'(1)-1|=\frac{1}{p}\ |P'_\lambda(1)+Q'(1)-1|=1.$$

 From Hensel's Lemma, there is an unique root of $p\ (Q_\lambda^*(z)-z)$ in $B_{r_0}(1)$, where $r_0= \frac{|Q(1)|}{p}$.
 We denote this root by $z_\lambda$ and observe that $z_\lambda $ is a fixed point of $Q_\lambda^*$.

 \n Now, we define the function

\begin{equation*}
\begin{array}{rccl}
h:& \Lambda& \longrightarrow& \ds \left\{z:|z-1|\leq \frac{|Q(1)|}{p}\right\} \\
& \lambda &  \longmapsto & \ \ z_\lambda
\end{array}
\end{equation*}
obtaining that

\begin{prop}{\hspace*{-6pt}\bf .}
The function $h$ is holomorphic in $\Lambda$.\label{res1}
\end{prop}
\n \D

\smallskip
\n Let $\{ h_{n} \}_{n\geq 0}$ be the sequence of functions
defined recursively as follows:
\begin{align*} &h_{0}(\lambda)=1\\
&h_{n}(\lambda)=
h_{n-1}(\lambda)-\frac{Q^{*}_{\lambda}(h_{n-1}(\lambda))}{(Q^{*}_{\lambda})'(h_{n-1}(\lambda))}\end{align*}

\n Then for all $n\in \N,\ h_{n}$ is a rational function without
poles in $\Lambda$.

\n As in the proof of Lemma \ref{Hensel}, we have
\begin{equation*}
\begin{array}{rcccl}
 | h_{n}(\lambda)-h(\lambda) |& = &\left| \ds\underset{i\geq
n}{\sum}\ds\frac{Q^{*}_{\lambda}(h_{i}(\lambda))}{(Q^{*}_{\lambda})'(h_{i}(\lambda))}\right|
&\leq& \ds\underset{i\geq n}{\max}
\left\{\left|\frac{Q^{*}_{\lambda}(h_{i}(\lambda))}{(Q^{*}_{\lambda})'(h_{i}(\lambda))}\right|
\right\}\\\\& =& \ds\underset{i\geq
n}{\max}{|Q^{*}_{\lambda}(h_{i}(\lambda) )|}&<& \rho ^{2^n}.
\end{array}
\end{equation*}

\n Hence $h_n$ converges to $h$ uniformly in $\Lambda$. Therefore,
$h\in \mathcal{H}(\Lambda).$ \hfill $\Box$

\bigskip

\bigskip
\n We may now introduce the affine map
$$A_\lambda(z)=z +h(\lambda)-1$$ we will work with the map

$$Q_\lambda(z)=A_\lambda^{-1}
(Q_\lambda^*(A_\lambda(z)))=P_\lambda(A_\lambda(z))+Q(A_\lambda(z))+1-h(\lambda).
$$ which is affinely conjugated to $Q_\lambda ^*$.
\bigskip
\n Notice that
\begin{equation*}
\begin{array}{rcl}
Q_\lambda(1)&=&P_\lambda(A_\lambda(1))+Q(A_\lambda(1))+1-h(\lambda)\\\\
&=&P_\lambda(h(\lambda))+Q(h(\lambda))+1-h(\lambda)\\\\
&=&1.
\end{array}
\end{equation*}

Moreover,
\begin{equation*}
\begin{array}{rcl}
|Q_\lambda'(1)|&=&|P'_\lambda(A_\lambda(1))+ Q'(A_\lambda(1))|\\\\
&=&|P'_\lambda(h(\lambda))+ Q'(h(\lambda))|\\\\
&=&p.
\end{array}
\end{equation*}

\n Thus, just as to the polynomials $P_\lambda$, $z=1$ is a
repelling fixed point of $Q_\lambda$ for all $\lambda \in
\Lambda$.

\bigskip
\n For the family $Q_\lambda$ we will obtain the following
theorem.
\bigskip
\begin{teo}{\hspace*{-6pt}\bf .}There is a dense set of parameter $\lambda \in \Lambda$ such that
the function $Q_{\lambda}$ has a wandering disc contained in the
filled Julia set, which is not attracted to an attracting
cycle.\label{teo que hago}
\end{teo}
\bigskip

\bigskip
\n {\bf Proof of Theorem \ref{teo pa estrella}.}

\n Recall that $Q_\lambda
(z)=A^{-1}_{\lambda}(Q^*_\lambda(A_\lambda(z)))$, i.e.
$Q_\lambda^*(z)= A_\lambda(Q_\lambda(A^{-1}_\lambda(z))).$

\n If $D$ is a ball, then $A_{\lambda} ^{-1}(D)$ and
$A_\lambda(D)$ are balls, and, obtaining directly that
$K(Q_\lambda,B)=A_\lambda^{-1}(K(Q_\lambda^*,A_\lambda(B)))$, it
is sufficient to show that if $D$ is a wandering disk for
$Q_\lambda$, then $A^{-1}(D)$ is a wandering disk for
$Q_\lambda^*.$

\n Suppose that $D$ is a wandering ball for $Q_\lambda$, i.e.
$Q_\lambda^m(D)\cap Q^n_\lambda(D) = \emptyset$ when $n\neq m$. It
follows that for $n\neq m$ we have that $A_\lambda^{-1}((Q_\lambda
*)^n(A_\lambda (D))) \cap A_\lambda^{-1}((Q_\lambda *)^m(A_\lambda
(D))) = \emptyset$, then $(Q_\lambda *)^n(A_\lambda
(D))\cap(Q_\lambda *)^m(A_\lambda (D))= \emptyset$. Therefore
$A_\lambda(D) $ is a wandering disk for $Q^*_\lambda$.\hfill$\Box$

\subsection{Properties of $\mathbf{Q_{\lambda}(}z\mathbf{)}$.}

\indent\indent The next proposition states a property of the
function $h$ that will be used several times.
\begin{prop}{\hspace*{-6pt}\bf .} If $\la{0},\la{1}
\in \La$, then $|h(\la{0})-h(\la{1})|\leq \rho
|\la{0}-\la{1}|$.\label{res2}
\end{prop}
\n \D

\n Let $\lambda _{0}, \lambda _{1}\in \Lambda$. Since
$$Q_{\lambda_0}^*(h(\lambda_0))-Q_{\lambda_1}^*(h(\lambda_1))=h(\lambda_0)-h(\lambda_1)$$
we have that
$$|Q_{\lambda_0}^*(h(\lambda_0))-Q_{\lambda_1}^*(h(\lambda_1))|\leq
\frac{|Q(1)|}{p}$$ and from
$$|Q(h(\la{0})) -Q(h(\la{1}))|<\ds \frac{\rho}{p} \
|h(\la{0})-h(\la{1})|,$$ it follows that
$$| P_{\lambda _{0}}(h(\la{0}))-P_{\la{1}}(h(\la{1}))|< \frac{\rho}{p} \
|h(\la{0})-h(\la{1})|.$$

In addition, from
$$| P_{\lambda _{0}}(h(\la{0}))-P_{\lambda _{0}}(h(\la{1}))|=p\ | h(\la{0})-h(\la{1}) |>
\frac{\rho}{p} \ |h(\la{0})-h(\la{1})|$$
$$| P_{\lambda
_{0}}(h(\la{1}))-P_{\lambda _{1}}(h(\la{1})) |=p\ | \lambda _{0}-
\lambda _{1} |\,
 |h(\la{1}) -1| $$ and by isosceles triangle principle, necessarily we have that $$p\ | \lambda
_{0}- \lambda _{1} | \,| h(\lambda_{1}) -1| =p\ |
h(\la{0})-h(\la{1}) |.$$

Finally from $|h(\lambda_1)-1|\leq \frac{\rho}{p}$, we have that
$| h(\la{0})-h(\la{1})|\leq \rho \,|\la{0}-\la{1}|.$\hfill $\Box$
\bigskip
\newpage
\begin{lema}{\hspace*{-6pt}\bf .} Let $z\in \mathbb{C}_p.$
\begin{itemize}
\item[ i)] If $\rho < |z |<1$, then $|Q_{\lambda}(z) |=p\ | z|^{p}
> |z|.$

\item[ii)] If $|z|\leq \rho$, then $|Q_{\lambda}(z)|\leq \rho.$

\item[iii)] If $|z-1|<1$, then $|Q_{\lambda}(z)-1|=p\ | z-1|.$
\item[iv)] If $1<|z|<\widehat{r}$, then $|Q_{\lambda}(z)|=p\
|z|^{p+1}$.
\end{itemize}\label{res3}
\end{lema}
\n \D

\n From Proposition \ref{res2}, for every $\lambda \in \Lambda$ we
have that $|h(\lambda)-1|\leq \frac{\rho}{p}< \rho$.
\begin{itemize}
\item[{\em i)}] Since $\rho<|z |<1$, we have that
$|A_{\lambda}(z)| =|z|$. In addition, $|A_{\lambda}(z)|^{p+1}<
|A_{\lambda}(z)|^{p}$. Hence $|P_{\lambda}(A_{\lambda}(z))|=p\
|A_{\lambda}(z) |^{p}>p\rho^{p}=\rho$. Furthermore,
$|Q(A_{\lambda}(z))|< \rho$ and $|1-h(\lambda)|<\rho$, therefore
 $$ |Q_{\lambda}(z) |=p\ | z |^{p}.$$

\item[{\em ii)}] Observe that $|A_{\lambda}(z)|\leq \rho$ since
$|z |\leq \rho $. It follows $$| P_{\lambda}(A_{\lambda}(z))|=p\
|A_{\lambda} (z)|^{p}\leq p\rho^{p}=\rho.$$ In addition, from
$|Q(A_{\lambda}(z))|< \rho$, $|1-h(\lambda)|<\rho$ and the strong
triangle inequality, we have that
$$|Q_{\lambda}(z)|\leq \rho.$$

\item[{\em iii)}]$$\begin{array}{rl}
|Q_{\lambda}(z)-1|&=|Q_{\lambda}(z)-Q_{\lambda}(1)|\\
&=|P_{\lambda}(A_{\lambda}(z))+ Q(A_{\lambda}(z))-
P_{\lambda}(A_{\lambda}(1))-Q(A_{\lambda}(1))|.\end{array}$$ Since
$|A_{\lambda}(z)-A_{\lambda}(1)|= |z-1|$, we have that
$|P_{\lambda}(A_{\lambda}(z))-
P_{\lambda}(A_{\lambda}(1))|=p\,|z-1|.$

Moreover,
$|Q_{\lambda}(A_{\lambda}(z))-Q_{\lambda}(A_{\lambda}(1))| \leq
\rho \,|z-1|$. Again, from the strong triangle inequality, we have
that
$$|Q_{\lambda}(z)-1|=p\ |z-1|.$$

\item[{\em iv)}] Since $|z|>1$, it follows that
$|P_{\lambda}(A_\lambda (z))|=p\ |z|^{p+1}$. Furthermore,
$|Q(A_\lambda (z))|< \rho$ and $|h(\lambda)-1|< \rho$, therefore
$$|Q_{\lambda}(z)|=p\ |z|^{p+1}.$$\hfill $\Box$
\end{itemize}

\n Recall that $B$ is the closed ball defined by $\{z \in \C_p
:|z| \leq \widehat{r}\}$, where $\widehat{r}$ is an element of
$|\C_p^*|$ chosen in the beginning of this section.
\begin{prop}{\hspace*{-6pt}\bf .}
 For each $\lambda \in \Lambda$, $(Q_{\lambda}, B)$ is an polynomial like map of degree $p+1$.\label{res4}
 \end{prop}
\n \D

\n Let $\lambda \in \Lambda$.  From the previous lemma we deduce
that $Q_{\lambda}(B)=\{z: |z|<p \ \widehat{r} ^{p+1}\}$, we will
prove that $(Q_{\lambda}, B)$ is of degree $p+1$.

\n Since $|P_\lambda(A_\lambda(z))-P_\lambda(z)|=p\
|h(\lambda)-1|<\rho$ we conclude that
$$Q_\lambda(z)-P_\lambda(z)=
P_\lambda(A_\lambda(z))-P_\lambda(z)+Q(A_\lambda(z))-h(\lambda)+1,$$
using that $\|Q_\lambda -P_\lambda \|< \rho$, the power series of
$Q_\lambda$ is
$$Q_{\lambda}(z)= a_0+a_1z+\ldots+\left (a_p+\frac{\lambda}{p}\right )z^p
 +\left(a_{p+1}+1-\frac{\lambda}{p}\right )z^{p+1}+a_{p+2}z^{p+2}\ldots$$
where  $\underset{i\geq 0}{\sup}\{ |a_i|z^{i}\}<\rho$. From
Theorem \ref{teorema de soluciones contar} it is possible to count
the solutions of $Q_{\lambda}(z)-w_0=0$.

\n If $p\leq|w_0|<pr^{p+1}$ and $f(z_0)=w_0$, then
$|z_0|=\left(\ds \frac{|w_0|}{p}\right )^{p+1}$, by Lemma
\ref{res3}{(\em iv)}. Therefore $w_0$ has $p+1$ pre-images in $B$.
Hence $(Q_{\lambda},B)$ is a polynomial like map of degree
$p+1$.\hfill $\Box$

\begin{prop}{\hspace*{-6pt}\bf .}
$K(Q_{\lambda},B)\subset B_{1}(0)\sqcup B_{1}(1).$\label{res5}
\end{prop}
\n \D

\n Suppose that $z\notin B_1(0) \sqcup B_1(1)$.

\n If $|z|>1$ then $|Q_{\lambda}(z)|= p|z|^{p+1}$. It follows that
there exists $n\in \mathbb{N}$, such that  $Q_{\lambda}^{n}(z)
\notin B.$

\n If $|z|=1$ and $|z-1|=1$, then
$|Q_{\lambda}(z)|=|P_{\lambda}(A_{\lambda}(z))+Q(A_{\lambda}(z))|=p$.
Hence $z\notin K(Q_{\lambda},B).$

\n Therefore $K(Q_{\lambda},B)\subset B_{1}(0)\sqcup
B_{1}(1).$\hfill $\Box$

\bigskip

\n This result allow us to define the itinerary of a point in
$K(Q_{\lambda},B)$. To simplify notation let $B_{0}=B_1(0)$ and
$B_1= B_1(1)$.

\n For any $z \in K(Q_{\lambda},B)$, the itinerary of $z$ for
$Q_{\lambda}$ is defined by
$$\theta _{0}\theta _{1}\ldots \theta _{n}\ldots  \in \{ 0,1\}^{\N \cup \{0\}} \text{\
where }Q_{\lambda}^{n}(z) \in B_{\theta _{n}} \text{\ for\ all \
}n \geq 0.$$

\begin{lema}{\hspace*{-6pt}\bf .}
Let $\lambda \in \Lambda $ and $D$ a ball contained in
$K(Q_{\lambda},B)$. Then:\label{res6}
\begin{itemize}
\item[1)]All points in $D$ have the same itinerary for
$Q_{\lambda}. $ \item[2)] If  the common itinerary of points in
$D$ is not pre-periodic, under the one side shift, then $D$ is a
wandering disc which is not attracted to an attracting periodic
point.
\end{itemize}
\end{lema}
\newpage \noindent \D

\begin{itemize}
\item[{\em1)}] We proceed by contradiction. Assume that there
exist $z_0\text{\ and\ }z_1\in D$ with different itineraries. Then
there exists $n_0 \in \mathbb{N}$ such that
$Q_{\lambda}^{n_{0}}(z_0)\in B_{0}$ and $
Q_{\lambda}^{n_{0}}(z_1)\in B_{1}$. Since $Q_{\lambda}^{n_{0}}(D)$
is a ball which has non-trivial intersection with $B_{0}$ and
$B_{1}$ we have that $\{z:| z | \leq 1 \} \subset
Q_{\lambda}^{n_{0}}(D)\subset K(Q_{\lambda},B)$, obtaining a
contradiction with Proposition \ref{res5}.

 \item[{\em2)}] Now, we suppose that $D$ is not a wandering disc, that is, there exist $n> m \geq 0$ such that $Q_{\lambda}^{n}(D) \cap Q_{\lambda}^{m}(D) \neq
 \emptyset$. Hence $Q_{\lambda}^{m}(z_0)\in Q_{\lambda}^{n}(D)$ for some $z_0\in D $. Therefore $Q_{\lambda}^{k}(z_0)\in
 Q_{\lambda}^{k+n-m}(D)$ for all $k$ in $ \mathbb{N}$.

 \n Since every point in $D$ has the same itinerary
 we conclude that the itinerary of the points in $D$ is pre-periodic with eventual period $n-m$.

 \n We must show that $D$ is not attracted to a periodic orbit. We suppose that there is an attracting periodic point $z_0$ and $s>0$ such that $B_s(z_0)$ is contained in the attracting
 basin of $z_0$, and $z_1\in D$ such that $z_1 \in
 B_s(z_0)$, from the first part of the proposition we have that
 every points of $B_s(z_0)$ have a common itinerary, and it is periodic.\hfill $\Box$

\end{itemize}

\n From the previous lemma we conclude that in order to prove
Theorem \ref{teo que hago} it is sufficient to find a wandering
disc in the filled Julia set of $(Q_\lambda,B)$ whose itinerary is
not pre-periodic, for a dense subset in $\Lambda$. Therefore, we
need to study the behavior of the points in $K(Q_\lambda,B)$ such
that its orbit visits both $B_0$ and $B_1$.

\n From Lemma \ref{res3} we know that the open ball $B_\rho(0)$ is
fixed under the action of $Q_\lambda$ and we have that a point
$x\in B$ has itinerary
$$\underset{n}{\underbrace{0\ldots0}}\,1\ldots$$ if and only if
$|x| =\rho_n$. Recall $\rho_0=1, \ p\rho^p_n=\rho_{n-1}$. This
crucial fact holds already for the family $P_\lambda$ \cite{RL2}.

\n The following lemma describe the local behavior of $Q_\lambda$
in the set $\{z:|z|=\rho_n\}$ and in $B_1$. \vspace*{-70pt}

\begin{lema}{\hspace*{-6pt}\bf .} $  \begin{array}{cl} \\*[59pt]
{\it 1}. & \text{Let } m\geq 1,|z_{0}|=| z_{1}|=\rho _{m}.\text{
If } | z_{0}-z_{1}| \leq S\text{, then}\\*[10pt] & \hspace*{3cm}
|Q_{\lambda}(z_{0})-Q_{\lambda}(z_{1})| \leq \rho _{m-1}
|z_{0}-z_{1}|.\\*[10pt] {\it 2}. & \text{ If } z_{0},z_{1} \in
B_{1}(1) \text{, then:}\\*[10pt] & \hspace*{3.4 cm
}|Q_{\lambda}(z_{0}) -Q_{\lambda}(z_{1})| =p\ |z_{0}-z_{1} |.
\end{array}$
\label{diferg}
\end{lema}
\n \D

\begin{itemize}
\item[{\em1)}]We observe that $|A_{\lambda}( z_{i})|=|z_{i}
|=\rho_{m}$, hence
 \\
$|Q(A_{\lambda}(z_{0}))- Q(A_{\lambda}(z_{1})) |< \rho \ |
A_{\lambda}(z_{0})-A_{\lambda}(z_{1})|=\rho \ |z_0-z_1|$.

Letting $\epsilon = A_{\lambda}(z_1)-A_{\lambda}(z_0)$ we have
\begin{equation*}
\begin{array}{rcl}
P_{\lambda}(A_{\lambda}(z_{0}))-P_{\lambda}(A_{\lambda}(z_{1}))
=&\ds \frac{\lambda}{p}\left (p\epsilon
A_{\lambda}(z_0)^{p-1}+\ldots
+p\epsilon^{p-1}A_{\lambda}(z_0)+\epsilon^p\right )\\\\
&+\left(1-\frac{\lambda}{p}\right )((p+1)\epsilon
A_{\lambda}(z_0)^p+\ldots +\epsilon ^{p+1}).
\end{array}
\end{equation*}
Moreover,
$|P_{\lambda}(A_{\lambda}(z_0))-P_{\lambda}(A_{\lambda}(z_1))|\leq
|\epsilon|\max \{A_{\lambda}(z_0)^{p-1},
p\,|\epsilon|^{p-1},\rho_{m-1}\}= |\epsilon|\ \rho_{m-1},$ since
 $|\epsilon|= |z_0-z_1|\leq S< \rho_m=|A_{\lambda}(z_0)|.$

Therefore $$|Q_{\lambda}(z_0)-Q_{\lambda}(z_1)|\leq \rho_{m-1}\
|z_0-z_1|.$$

 \item[{\em2)}] We note that if $y \in B_{1}$, then
$A_{\lambda}(y) \in B_{1}$. Now
\begin{equation*}
\begin{array}{rcl}
P_{\lambda}(A_{\lambda}(z_{0}))-P_{\lambda}(A_{\lambda}(z_{1}))
=&\ds \frac{\lambda}{p}\left (p\,\epsilon
A_{\lambda}(z_0)^{p-1}+\ldots
+p\, \epsilon^{p-1}A_{\lambda}(z_0)+\epsilon^p\right )\\\\
&+\left(1-\frac{\lambda}{p}\right )((p+1)\,\epsilon \,
A_{\lambda}(z_0)^p+\ldots +\epsilon ^{p+1}),
\end{array}
\end{equation*}
where $\epsilon = z_0 -z_1 $, thus
$|P_{\lambda}(A_{\lambda}(z_{0}))-P_{\lambda}(A_{\lambda}(z_{1})|=p\
| z_{0}-z_{1}|$ since $|\epsilon|<1$. Furthermore we have that
$$|Q(A_{\lambda}(z_{0}))-Q(A_{\lambda}(z_{1})|\leq \rho \
|z_{0}-z_{1} | <p\ |z_{0}-z_{1} |.$$ Therefore
$$|Q_{\lambda}(z_{0})-Q_{\lambda}(z_{1})|= p\ |
z_{0}-z_{1}|.$$\hfill $\Box$
\end{itemize}
\bigskip

\n The following lemma gives a sufficient condition for the
existence of a wandering disc in $K(Q_\lambda,B)$.

\begin{lema}{\hspace*{-6pt}\bf .}
Let $\{ m_{i}\} _{i\geq 0},\{ M_{i}\} _{i\geq 0}$ be sequences of
positive integers such that if $i\geq0$
$$\rho_{m_{i}-1}\cdot ...\cdot \rho_{1}\cdot p^{M_{i}}\leq1.$$

\n If for $\lambda_{0} \in \Lambda  $ there exists $z_0\in
K(Q_{\lambda_{0}})$ with itinerary
$$\underset{m_0}{\underbrace{0\ldots0}}\underset{M_0}{\underbrace{1\ldots1}}\underset{m_
1}{\underbrace{0\ldots0}}
\underset{M_1}{\underbrace{1\ldots1}}\ldots\underset{m_i}{\underbrace{0\ldots0}}\underset{M_i}{\underbrace{1\ldots1}}\ldots,$$
then the closed ball $D=\{|z-z_0|\leq S \}$ is contained in
$K(Q_{\lambda_{0}})$.\label{res8}

\n If we add the hypothesis  $ \lim M_i= \infty$, then by Lemma
\ref{res6} we have that $D$ is a wandering disc contained in
$K(Q_{\lambda},B)$ which is not attracted to an attracting cycle.
\end{lema}
\n \D

 \n We now define the sequence $\{ N_{i}\}_{i\geq 0}$
recursively:
$$N_{0}=0$$
$$N_{i}=N_{i-1}+m_{i-1}+M_{i-1}.$$

\n We will prove inductively that $\diam (Q_{\lambda
_{0}}^{N_{i}}(D))\leq S$. For $i=0$ the claim is true because the
definition of $D$. \n Suppose that $\diam (Q_{\lambda
_{0}}^{N_{i}}(D))\leq S$, since $Q_{\lambda
_{0}}^{N_{i}+j}(z_0)\in B_{0}$ for all $ j$ in $\N$ \linebreak
with $0\leq j\leq m_i $ and $Q_{\lambda _{0}}^{N_i+m_{i}}(z_0) \in
B_{1}$, it follows that $|Q_{\lambda _{0}}^{N_{i}}(z_0)| =\rho
_{m_{i}}$. Therefore $|Q_{\lambda _{0}}^{N_{i}}(y)| =\rho
_{m_{i}}$ and $Q_{\lambda _{0}}^{N_i+m_{i}}(y) \in B_{1} $ for all
$ y \in D.$

\n  From the first statement of the previous lemma we have that
$$\diam(Q_{\lambda _{0}}^{N_{i}+m_{i}}(D))\leq \rho _{m_{i-1}}\cdot ...\cdot \rho _{1}
\diam(Q_{\lambda _{0}}^{N_{i}}(D))\leq p^{-M_{i}}S.$$ \n Now,
using the second statement of the same lemma, we obtain
$$\diam(Q_{\lambda _{0}}^{N_{i+1}}(D))\leq p^{M_{i}}\diam(Q_{\lambda
_{0}}^{N_{i}+m_{i}}(D))\leq S.$$ \n Therefore $D\subset
K(Q_{\lambda _{0}},B)$. \hfill $\Box$

\newpage
\subsection{Parameter selection.}

\indent\indent In this section we will prove results that describe
the behavior of the iterates $Q_{\lambda}^{n}(z) $, not just as a
function of $z$ but also as a function of $\lambda$.
\begin{lema}{\hspace*{-6pt}\bf .}
Let $\lambda_0,\lambda_1 \in \Lambda$. \label{res9}
\begin{itemize}
\item[1)] If $z\in B_0$ then
$|P_{\lambda_0}(z)-P_{\lambda_1}(z)|=p\
|z|^p|\lambda_0-\lambda_1|.$ \item[2)] If $z\in B_1$ then
$|P_{\lambda_0}(z)-P_{\lambda_1}(z)|=p\ |\lambda_0-\lambda_1|\
|z-1|.$
\end{itemize}
\end{lema}
\n \D

\begin{itemize}
\item[{\em1)}]
\begin{align*}
&|P_{\lambda_0}(z)-P_{\lambda_1}(z)| =\left |z^p\left (\ds
\frac{\lambda_0-\lambda_1}{p}\right )-z^{p+1}\left
(\ds\frac{\lambda_0-\lambda_1}{p}\right )\right | =p\
|\lambda_0-\lambda_1|\ |z-1|.
\end{align*}

\item[{\em 2)}]
\begin{align*}
&|P_{\lambda_0}(z)-P_{\lambda_1}(z)| =\left|z^p\left (\ds
\frac{\lambda_0-\lambda_1}{p}\right )-z^{p+1}\left (\ds
\frac{\lambda_0-\lambda_1}{p}\right )\right | =p\ |x|^p\
|\lambda_0-\lambda_1|.
\end{align*}
 \hfill $\Box$
\end{itemize}
\begin{lema}{\hspace*{-6pt}\bf .}
Let $M\in \mathbb{N}$ and $x_{0},x_{1}\in \{ x:|x-1 | \leq
p^{-M}\}$. If the parameters $\lambda _{0}, \lambda _{1}\in
\Lambda$\label{res10} are such that $| \lambda _{0}- \lambda _{1}|
= | x_{0}-x_{1}|$, then
$$| Q_{\lambda_{0}}^{M}(x_{0})-Q_{\lambda_{1}}^{M}(x_{1})| =p^{M} |\lambda _{0}- \lambda 
_{1} |.$$
\end{lema}
\n \D \setcounter{equation}{0} \ First we prove the lemma for
$M=1$.

\n By the second part of Lemma \ref{diferg} we have

\begin{equation}| Q_{\lambda _{0}}(x_{0})-Q_{\lambda
_{0}}(x_{1})  |= p\  | x_{0}-x_{1}|=p \ |
\lambda_{0}-\lambda_{1}|.\label{c1}
\end{equation}
\n Furthermore by Lemma \ref{comp de P} and Lemma \ref{res2} we
obtain
\begin{equation}| P_{\lambda _{0}}(A_{\lambda_0}(x_{1}))-P_{\lambda
_{0}}(A_{\lambda_1}(x_{1}))| =p\  | h(\lambda _{0})-h(\lambda
_{1})| < p \ |\lambda _{0}-\lambda_{1}|.\label{c2}
\end{equation}

\n By Lemma \ref{res9}, we conclude that
\begin{equation}
| P_{\lambda _{0}}(A_{\lambda_1}(x_{1}))-P_{\lambda
_{1}}(A_{\lambda_1}(x_{1}))|=p\ |\lambda _{0}-\lambda_{1}|
\,\,|A_{\lambda_1}(x_{1})-1 |< p \ |\lambda
_{0}-\lambda_{1}|,\label{c3}
\end{equation}
\n and from equations (\ref{c2}) and (\ref{c3}) \begin{equation} |
P_{\lambda _{0}}(A_{\lambda _{0}}(x_{1}))-P_{\lambda
_{1}}(A_{\lambda _{1}}(x_{1}))|<p\  |\lambda _{0}-\lambda _{1}|.
\label{c4}
\end{equation}

\n Moreover,
\begin{equation}
 |Q(A_{\lambda_0}(x_{1}))-
Q(A_{\lambda_1}(x_{1}))|< \rho \ | h(\lambda _{0})-h(\lambda _{1})
| <p\  |\lambda _{0}-\lambda _{1}|. \label{c5}\end{equation}

\n From (\ref{c4}) and (\ref{c5}) we have that
\begin{equation}
| Q_{\lambda _{0}}(x_{1})-Q_{\lambda _{1}}(x_{1})  |< p \ |
\lambda_{0}-\lambda_{1}|. \label{c6}
\end{equation}

\n Finally, from (\ref{c1}) and (\ref{c6}) we obtain that $|
Q_{\lambda_{0}}(x_{0})-Q_{\lambda _{1}}(x_{1}) |= p \ |
\lambda_{0}- \lambda_{1}|. $

\n Now let us prove that the proposition is true for $M+1$. By the
inductive hypothesis and the third statement of Lemma \ref{res3}
we have that
$Q_{\lambda_{0}}^{M}(x_{0}),Q_{\lambda_{1}}^{M}(x_{1})$ belong to
$B_{1}$ and using Lemma \ref{diferg} with
$z_0=Q_{\lambda_0}^M(x_0) $ and $z_1=Q_{\lambda_1}^M(x_1)$ we
obtain that
\begin{equation}
|
Q_{\lambda_{1}}(Q_{\lambda_{0}}^{M}(x_{0}))-Q_{\lambda_{1}}(Q_{\lambda_{1}}^{M}(x_{1}))|=
p^{M+1}| \lambda _{0}- \lambda _{1}|.\label{c7}
\end{equation}

\n Furthermore
\begin{equation}
|P_{\lambda_{0}}(A_{\lambda_0}(Q_{\lambda_{0}}^{M}(x_{0})))-P_{\lambda_{0}}(A_{\lambda_1
}(Q_{\lambda_{0}}^{M}(x_{0})))|=p\
\label{c8}|h(\lambda_{0})-h(\lambda_{1})| <p \
|\lambda_{0}-\lambda_{1} |,
\end{equation}
just as before, from the previous lemma
\begin{equation}
|P_{\lambda_{0}}(A_{\lambda_1}(Q_{\lambda_{0}}^{M}(x_{0})))-P_{\lambda_{1}}(A_{\lambda_1
}(Q_{\lambda_{0}}^{M}(x_{0})))|<p\
\label{c9}|\lambda_{0}-\lambda_{1}|
\end{equation}
\n and

\begin{equation}
| Q(A_{\lambda_0}\label{c10}(Q_{\lambda_{0}}^{M}(x_{0})))-
Q(A_{\lambda_1}(Q_{\lambda_{0}}^{M}(x_{0}))) |< \rho \
|h(\lambda_{0})-h(\lambda_{1})|<p \ |\lambda_{0}-\lambda_{1}|.
\end{equation}
The strong triangle principle applied to (\ref{c8}), (\ref{c9})
and (\ref{c10}) gives us
\begin{equation}
|Q_{\lambda_0}(Q_{\lambda_0}^M(x_0))-Q_{\lambda_1}(Q_{\lambda_0}^M(x_0))|<p\
|\lambda_0-\lambda_1|, \label{c11}
\end{equation}
and from (\ref{c7}) and (\ref{c11}) we conclude that
$$| Q_{\lambda_{0}}^{M+1}(x_{0})-Q_{\lambda_{1}}^{M+1}(x_{1})|
=p^{M+1} |\lambda _{0}- \lambda _{1} |.$$ \hfill $\Box$

\begin{lema}{\hspace*{-6pt}\bf .}
Let $m \in \mathbb{N}$ and let $x_{0},x_{1}$ be such that
$|x_0|=|x_1|=\rho _{m}$ and $| x_{0}-x_{1}|\leq S$. If $\lambda
_{0}, \lambda _{1}\in \Lambda$ are such that $$\rho _{m-1} \cdot
\ldots \cdot \rho_1| x_{0}-x_{1}| < | \lambda _{0}- \lambda _{1}|
\leq S,$$ then
$$| Q_{\lambda_{0}}^{m}(x_{0})-Q_{\lambda_{1}}^{m}(x_{1})|= | \lambda _{0}- \lambda
_{1}|.$$\label{res11}
\end{lema}
\n \D

We start by inductively proving that if $1\leq i\leq m $, then
$$| Q_{\lambda_{0}}^{i}(x_{0})- Q_{\lambda_{1}}^{i}(x_{1}) | \leq \max \{  \rho_{m-i}|
\lambda_{0}-\lambda_{1}|, \rho_{m-1} \cdot \ldots \cdot \rho _{m-i}| x_{0}-x_{1}| \}.$$

From the first part of Lemma \ref{diferg} we have
\begin{equation}| Q_{\lambda_{0}}(x_{0})-Q_{\lambda_{0}}(x_{1})|
\leq \rho_{m-1}| x_{0}-x_{1}|. \label{c12}
\end{equation}

Since $|A_{\lambda_0}( x_{1})|=|A_{\lambda_1}( x_{1})|= \rho
_{m}$, and $\rho _{m-1} \cdot \ldots \cdot \rho_{1}|
h(\lambda_{0})-h(\lambda_{1})| <| \lambda_{0}-\lambda_{1}| \leq
S$, we obtain
\begin{equation}
|P_{\lambda_0}(A_{\lambda_0}(x_1))-P_{\lambda_1}(A_{\lambda_1}(x_1))|
\leq \max \{\rho_{m-1}|h(\lambda_0)-h(\lambda_1)|,
\rho_{m-1}|\lambda_0-\lambda_1|\}, \label{c13}
\end{equation}
by the equation (\ref{eq6}) of Lemma \ref{lema en P inductivo}.

\n Moreover
\begin{equation}
|Q(A_{\lambda_0}(x_{1}))-Q(A_{\lambda_1}(x_{1}))|<\rho \,|
h(\lambda_{0})-h(\lambda_{1})|<\rho\,|\lambda_0-\lambda_1|.
\label{c14}
\end{equation}
From inequalities (\ref{c13}) and (\ref{c14}), together with
Proposition \ref{res2} we have
\begin{equation}
\label{c16}| Q_{\lambda_{0}}(x_{1})-Q_{\lambda_{1}}(x_{1})|\leq
\rho_{m-1}|\lambda_0 -\lambda _1|.
\end{equation}
Now inequalities (\ref{c12}) and (\ref{c16}) give us
$$| Q_{\lambda_{0}}(x_{0})-Q_{\lambda_{1}}(x_{1})|\leq \max \{\rho_{m-1}\, |x_0-x_1|,
\rho_{m-1}\, |\lambda_0-\lambda_1|\}.$$

\n Therefore which one is true for $i=1$.

\n Now suppose $| Q_{\lambda_{0}}^{i}(x_{0})-
Q_{\lambda_{1}}^{i}(x_{1}) | \leq \max \{  \rho_{m-i}|
\lambda_{0}-\lambda_{1}|, \rho_{m-1} \cdot \ldots \cdot \rho
_{m-i}| x_{0}-x_{1}| \}$.

\n Notice that
$|Q_{\lambda_{0}}^{i}(x_{0})|=|Q_{\lambda_{1}}^{i}(x_{1})|=\rho
_{m-i}$, therefore, using  Lemma \ref{diferg} with
$z_0=Q_{\lambda_{0}}^{i}(x_{0})$ and
$z_1=Q_{\lambda_{1}}^{i}(x_{1})$, we obtain that
\begin{equation}
 |
Q_{\lambda _{1}}(Q_{\lambda_{0}}^{i}(x_{0}))-Q_{\lambda
_{1}}(Q_{\lambda_{1}}^{i}(x_{1}))|\leq \rho _{m-(i+1)} |
Q_{\lambda_{0}}^{i}(x_{0})- Q_{\lambda_{1}}^{i}(x_{1})
|.\label{c17}
\end{equation}

\n From Lemma \ref{res9} with
$x=A_{\lambda_0}(Q_{\lambda_{0}}^{i}(x_{0}))$ we have that
\begin{equation} |
P_{\lambda_{0}}(A_{\lambda_0}(Q_{\lambda_{0}}^{i}(x_{0})))-P_{\lambda_{1}}(A_{\lambda_0}
(Q_{\lambda_{0}}^{i}(x_{0})) )|=p\,|
\lambda_{0}-\lambda_{1}|\,\rho_{m-i}^{p}= \rho_{m-(i+1)}|
\lambda_{0}-\lambda_{1}|.\label{c18}
\end{equation}
Using the first part of Lemma \ref{comp de P} we obtain that
\begin{equation*}
 |P_{\lambda_{1}}(A_{\lambda_0}(Q_{\lambda_{0}}^{i}(x_{0})))-P_{\lambda_{1}}(A_{\lambda_1}
(Q_{\lambda_{1}}^{i}(x_{1})) )|\\ \leq \rho _{m-(i+1)}|
Q_{\lambda _{0}}^{i}(x_{0})-Q_{\lambda
_{1}}^{i}(x_{1})+h(\lambda_{0})-h(\lambda_{1})|,
\end{equation*}
\begin{equation*}
|Q(A_{\lambda_0}(Q_{\lambda_0}^{i}(x_0)))-Q(A_{\lambda_1}(Q_{\lambda_0}^{i}(x_0)))|<\rho
\,|z_{\lambda_0}-z_{\lambda_1}|<\rho\,|\lambda_0-\lambda_1|.
\end{equation*}

\n From the inequalities above and Proposition \ref{res2} we have
that
 $$| Q_{\lambda_{0}}^{i+1}(x_{0})- Q_{\lambda_{1}}^{i+1}(x_{1})
| \leq \max \{  \rho_{m-(i+1)}| \lambda_{0}-\lambda_{1}|,\,
\rho_{m-1} \cdot \ldots \cdot \rho _{m-(i+1)}| x_{0}-x_{1}| \}.$$

\n Notice that in the inductive step for $i=m-1$, we have that
$$|Q_{\lambda_{0}}^{m-1}(x_{0})- Q_{\lambda_{1}}^{m-1}(x_{1}) | \leq
\max \{  \rho_{1}| \lambda_{0}-\lambda_{1}|, \rho_{m-1} \cdot
\ldots \cdot \rho _{1}| x_{0}-x_{1}| \}.$$ From the above
inequality we obtain
 $$|Q_{\lambda_{1}}(Q_{\lambda_{0}}^{m-1}(x_{0}))-Q_{\lambda
_{1}}(Q_{\lambda_{1}}^{m-1}(x_{1}))|\leq \rho _{1} |
Q_{\lambda_{0}}^{m-1}(x_{0})- Q_{\lambda_{1}}^{m-1}(x_{1}) |< |
\lambda_{0}-\lambda_{1}|$$ and from lemmas \ref{res9} and
\ref{comp de P} we have the following inequalities
\begin{equation*}
\begin{array}{rcl}
 |
P_{\lambda_{1}}(A_{\lambda_0}(Q_{\lambda_{0}}^{m-1}(x_{0})))-P_{\lambda_{1}}(A_{\lambda_
1}(Q_{\lambda_{1}}^{m-1}(x_{1})) )|&\leq &|  Q_{\lambda
_{0}}^{m-1}(x_{0})-Q_{\lambda
_{1}}^{m-1}(x_{1})+h(\lambda_{0})-h(\lambda_{1})|\\\\
&<& | \lambda_{0}-\lambda_ {1}|.
\end{array}
\end{equation*}

\begin{equation*}|
P_{\lambda_{0}}(A_{\lambda_0}(Q_{\lambda_{0}}^{m-1}(x_{0})))-P_{\lambda_{1}}(A_{\lambda_
0}(Q_{\lambda_{0}}^{m-1}(x_{0})) )|=p\,|
\lambda_{0}-\lambda_{1}|\,\rho_{1}^{p}= | \lambda_{0}-\lambda_{1}|
\end{equation*}

Moreover
\begin{equation*}
| Q(A_{\lambda_0}(Q_{\lambda_{0}}^{m-1}(x_{0})))-
Q(A_{\lambda_1}(Q_{\lambda_{0}}^{m-1}(x_{0})))|< \rho |
h(\lambda_{0})-h(\lambda_{1})|<| \lambda _{0}-\lambda_{1}|
\end{equation*}
and using the four previous inequalities and Proposition
\ref{res9} we have
$$| Q_{\lambda_{0}}^{m}(x_{0})-Q_{\lambda_{1}}^{m}(x_{1})|= |
\lambda _{0}- \lambda _{1}|. $$ \hfill$\Box$
\bigskip
\begin{prop}{\hspace*{-6pt}\bf .}
Let $\lambda \in \Lambda$ and consider $x\in K(Q_{\lambda},B)$
with itinerary  $$\theta _{0}\,
\theta_{1}\ldots\theta_{n-1}\,1\,1\,1\ldots,$$ for $Q_\lambda$,
i.e. $Q^{n}_{\lambda}(x)=1$, for some $n\geq 1$.

\n Suppose that there exists $\epsilon \in (0,1)$ such that for
all $\lambda_{0},\lambda_{1}$ in $\{\omega:| \omega -\lambda|\leq
\epsilon \}$, is true that
$|Q^{n}_{\lambda_{0}}(x)-Q^{n}_{\lambda_{1}}(x)|=|
\lambda_{0}-\lambda_{1}|.$

Let $M,m \in \mathbb{N}$ be such that $p^{-M}\leq \epsilon $ and
$$p^{M} \cdot \rho _{m-1}\cdot\ldots\cdot \rho_{1}< 1.$$
Then there exists $\lambda' \in \Lambda$ with $|\lambda -\lambda'|
\leq p^{-M}$ such that $x$ has itinerary
$$\theta _{0}\, \theta_{1}\ldots\theta_{n-1}\,\underset{M}{\underbrace{1\ldots
1}}\,\underset{m}{\underbrace{0\ldots 0}}\,1\,1\ldots\text{ for }
Q_{\lambda '}$$ and such that for all pairs of elements
$\lambda_{0},\lambda_{1} $ in $\{\omega:| \omega -\lambda' |\leq
Sp^{-M} \}$, we have that
$$ | Q_{\lambda_{0}}^{n+m+M}(x)-Q_{\lambda_{1}}^{n+m+M}(x)|=|
\lambda_{0}-\lambda_{1}|.$$\label{paso inductivo}
\end{prop}
\D \  Let $\phi: \Lambda\longrightarrow \C_p$ be the function
defined by $\phi(w)= Q_{w}^{n+M}(x)$. By Proposition \ref{res1} we
have that $\phi$ is holomorphic in $\Lambda$. Furthermore, by
hypothesis, if $\lambda_{0},\lambda_{1}\in B_{p^{-M}}(\lambda)$,
then $|Q_{\lambda_{0}}(x)-Q_{\lambda_{1}}(x)|=| \lambda_{0}-
\lambda_{1}|\leq p^{-M }.$

\n Applying  Lemma \ref{res10} to $z_0=Q_{\lambda_{0}}^n(x)$ and
$z_1=Q_{\lambda_{1}}^n(x)$, we have that
\begin{equation}| \phi(\lambda_{0}) -\phi(\lambda_{1}) |= p^{M}
|\lambda_{0}- \lambda_{1}|.\label{c20}\end{equation}

Then, by Corollary \ref{esfera}, we have that $\phi (\{ w:|
w-\lambda|= p^{-M}\})= \{ w:| w-1|= 1\}$. Therefore, there exists
$w_{0}\in \Lambda$ such that $\phi(w_{0})=0$.

\n If $w\in \{z: |z-w_0|< p^{-M}\}$, the itinerary of $x$ for
$Q_w$ is
$$\theta _0 \, \theta_{1}\ldots\theta_{n-1}\,\underset{M}{\underbrace{1\ldots
1}}\,0\ldots $$ this is direct consequence of (\ref{c20}).

\n Now let us consider the function $\psi:\Lambda\longrightarrow
\C_p$ defined by $\psi(w)=Q_{w}^{n+M+m}(x)$.

\n Since $Q_\lambda$ leaves $B_\rho(0)$ fixed, we obtain that
$|\psi(w_0)|<\rho$. Now, by (\ref{c20}) we have that
$|\phi(w)|=\rho_m$ if $w$ is such that $|w-w_0|=\rho_m p^{-M}$,
hence $|\psi(w)|=1$. Using again Corollary \ref{esfera}, we
observe that
$$\psi(\{ w:| w-w_{0}|= p^{-M}\rho_{m}\})= \{ w:| w|= 1\}.$$

\n Therefore, there exists $\lambda'$ such that
$\psi(\lambda')=1$.

\n Thus, the itinerary of $x$ for $Q_{\lambda'}$ is
$$\theta _0 \, \theta_{1}\ldots\theta_{n-1}\,\underset{M}{\underbrace{1\ldots
1}}\underset{m}{\underbrace{0\ldots 0}}\,1\,1\ldots$$

\n Notice that $|\phi(\lambda')|=\rho_m$. If $\lambda_0, \lambda_1
$ belong to $\{w:|w-\lambda'|\leq S\,p^{-M}\}$, then the points
$z_0=\phi(\lambda_0)$ y $z_1=\phi(\lambda_1)$ are such that
$|z_0|=|z_1|=\rho_m$ and $|z_0-z_1|=p^M |\lambda_0-\lambda_1|\leq
S$, by (\ref{c20}). Moreover, by hypothesis, we have
$$\rho_{m-1}\cdots \rho_1 |z_0-z_1|=\rho_{m-1}\cdots\rho_1 p^M
|\lambda_0-\lambda_1|< |\lambda_0-\lambda_1|.$$

\n Now, applying Lemma \ref{res11} we obtain that
$$|Q_{\lambda_0}^{n+M+m}(x)-Q_{\lambda_1}^{n+M+m}(x)|=|\lambda_0-\lambda_1|.$$
\hfill$\Box$

\noindent {\bf Proof of Theorem \ref{teo que hago}.}

\n We define the sequence $\{M_i\}_{i \in \N}$ recursively. Choose
$M_0 \in \N$ such that $p^{-M_0} \leq S$, and suppose that $M_i$
is already defined. Now choose $M_{i+1} \in \N$ satisfying
$p^{M_{i+1}-M_i}\leq S$.

\n Furthermore, we define $\{m_i\}_{i\in \N}$ such that for each
$i \in \N$,
$$\rho _{m_i -1}\cdot\ldots\cdot\rho_1 \cdot p^{M_i}\leq 1.$$

\n For an arbitrary $\lambda_0 \in \Lambda$ there exists $x \in
B_0$ such that $Q_{\lambda_0}^{m_0}(x)=1$, i.e. its itinerary for
$Q_{\lambda_0}$ is
$$\underset{m_0}{\underbrace{0\ldots0}}\,1\,1\,1\ldots$$

\n By Lemma \ref{res11}, for $\lambda \in \Lambda$ with $|\lambda-
\lambda_0|\leq S$, we have
$$|Q_{\lambda}^{m_0}(x)-Q_{\lambda_0}^{m_0}(x)|=|\lambda-\lambda_0|.$$

Since $\rho_{m_{i+1}-1}\cdot\ldots\cdot\rho_1\cdot p^{M_{i+1}}\leq
1$, we have that
$$p^{M_i}\cdot\rho_{m_{i+1}-1}\cdot\ldots\cdot\rho_1\leq S <1.$$

\n Therefore for $\lambda=\lambda_0,n=m_0, m=m_1$ and $\epsilon
=p^{-M_0}$ the hypothesis of Proposition \ref{paso inductivo}
hold. Hence  we may consider $\lambda_1$ with
$|\lambda_0-\lambda_1|\leq p^{-M_0}$ and such that the itinerary
of $x$ for $Q_{\lambda_1}$ is
$$\underset{m_0}{\underbrace{0\ldots 0}}\, \underset{M_0}{\underbrace{1\ldots 1}}\,
\underset{m_1}{\underbrace{0\ldots 0}}\, 1\, 1\ldots.$$

\n In view of the second part of Proposition \ref{paso inductivo},
for all pairs of elements $\omega_0, \omega_1$ in
$\{\omega:|\omega - \lambda_1|\leq Sp^{-M_0} \}$ we have that
$|Q_{\omega_0}^{m_0+M_0+m_1}(x)-Q_{\omega_1}^{m_0+M_0+m_1}(x)|=|\omega_0-\omega_1|$,
then we can use this proposition recursively. For the i-th step we
consider $n=n_0+M_0+\ldots+m_{i-1}+M_{i-1}+m_{i},\, \lambda =
\lambda_{i-1}$ y $\epsilon_i =Sp^{-M_{i-1}}$, obtaining
 $\lambda_i \in \Lambda$ with $|\lambda_i-\lambda_{i-1}|\leq
Sp^{-M_i}$ and such that the itinerary of $x$ for $Q_{\lambda_i}$
is
$$\underset{m_0}{\underbrace{0\ldots0}}\,\underset{M_0}{\underbrace{1\ldots1}}\ldots
\underset{M_{i-1}}{\underbrace{1\ldots1}}\,
\underset{m_i}{\underbrace{0\ldots0}}\,1\,1\,1\ldots$$

\n By definition $\underset{i\rightarrow \infty}{\lim}M_i=
\infty$, and since $|\lambda_{i+1}-\lambda_i|\leq Sp^{-M_i}$,
 $\{\lambda_i\}_{i\in \N}$ is a Cauchy sequence. If we call its limit $\lambda$ we have that $|\lambda-\lambda_0|\leq
p^{-M_0}$ and the itinerary of $x$ for $Q_{\lambda}$ is
$$\underset{m_0}{\underbrace{0\ldots0}}\,\underset{M_0}{\underbrace{1\ldots1}}\ldots
\underset{M_{i-1}}{\underbrace{1\ldots1}}\,
\underset{m_i}{\underbrace{0\ldots0}}\,\underset{M_{i}}{\underbrace{1\ldots1}}\ldots.$$

Moreover the sequences $\{M_i\}_{i\in\N},\{m_i\}_{i\in\N}$ satisfy
the hypothesis of Lemma \ref{res8}, therefore $Q_\lambda$ has a
wandering disc contained in $K(Q_{\lambda},B)$, which is not
attracted to an attracting cycle.

\n Finally, recall that $\lambda_0\in \Lambda$ and $M_0 \in \N$
were chosen arbitrarily and since $|\lambda-\lambda_0|\leq
p^{-M_0}$ we have that for a dense set of parameters $\lambda \in
\Lambda$ the function $Q_\lambda$ has a wandering disc in
$K(Q_\lambda, B)$ which is not attracted to an attracting
cycle.\hfill$\Box$

\end{document}